\documentclass[12pt,a4paper]{article}
\usepackage{amssymb,amsmath}
\usepackage{amsfonts}
\usepackage{centernot}
\usepackage{mathtools}
\usepackage{stmaryrd}
\usepackage{a4}
\newcommand{\para}{\par\vspace{.25cm}}

\newtheorem{prop}{Proposition}
\newtheorem{theorem}{Theorem}

\newtheorem{cor}{Corollary}

\usepackage{multirow}

\begin{document}
\baselineskip 18pt

\title{\bf Integral group rings with all central units trivial: solvable groups}
\author{Sugandha Maheshwary {\footnote {Research supported by SERB, India (PDF/2016/000731).}}
\\ {\em \small Indian Institute of Science Education and Research, Mohali,}\\
{\em \small Sector 81, Mohali (Punjab)-140306, India.}
\\{\em \small email: sugandha@iisermohali.ac.in}}
\date{}
{\maketitle}

\begin{abstract}\noindent {}
The object of this paper is to examine finite solvable groups whose integral group rings have only trivial central units.
\end{abstract}\vspace{.25cm}
 \vspace{.25cm} 

\section{Introduction}
Let $G$ be a finite group. It is easy to see that the group $\mathcal{Z}(\mathcal{U}(\mathbb{Z}[G]))$ of central units of the integral group ring $\mathbb{Z}[G]$ contains the \emph{trivial central units} $\pm g$, $g \in \mathcal{Z}(G)$, the centre of $G$. Following \cite{BMP}, we say that $G$ is a group with the {\sf cut}-property, if $\mathcal{Z}(\mathcal{U}(\mathbb{Z}[G]))$ consists of only the trivial central units. Nilpotent groups and metacyclic groups  with the {\sf cut}-property have been investigated and a complete classification of finite metacyclic groups with this property has been given in \cite{BMP}.\linebreak Continuing this investigation, and drawing on the work in \cite{CD}, we provide, in \linebreak Section 2, a precise characterization for a solvable group $G$ to have the {\sf cut}-property in case (i) the order of $G$ is odd, or (ii) the order of  $G$ is even and every element of $G$ has prime-power order (Theorems \ref{T5} \& \ref{T7}). We also complete the classification of finite nilpotent groups with the {\sf cut}-property (Theorem \ref{T4}).

\section{Solvable groups with the {\sf cut}-property}

Throughout the paper, $G$ denotes a finite group and for elements $x,y \in G$, $x\sim y$, denotes ``\emph{x is conjugate to y}." We recall (\cite{sehg} \& \cite{DMS}, Lemma~2) that $G$ has the \linebreak{\sf cut}-property, if, and only if, for every $x\in G$ and for every integer $j$ relatively prime to the order $o(x)$ of $x$,  \begin{equation}\label{Eq4a}
    x^{j}\sim x ~\mathrm{or} ~x^{-1}.
\end{equation}
 \begin{theorem}\label{T5}An odd-order group $G$ has the {\sf cut}-property if, and only if, every \linebreak element $x \in G$ satisfies {\rm(i)} $x^{5}\sim x^{-1}$, and {\rm(ii)} $o(x)$ is either $7$ or a power of~$3$.
 \end{theorem}\textbf{Proof.} Let $G$ be an odd-order group with the {\sf cut}-property, so that $G$ satisfies~(\ref{Eq4a}), i.e., in the terminology of Chillag-Dolfi \cite{CD}, $G$ is an \emph{inverse semi-rational group}. Consequently, by (\cite{CD}, Theorem 3) together with an observation in its proof that the order of every element in such a group is a prime power and (\ref{Eq4a}), we have that for every element $x\in G$,

 \begin{equation}\label{Eq4}
    o(x)=3^{a} ~{\mathrm or}~ 7^{a}, ~a\geq 0
\end{equation}

\noindent and
 \begin{equation}\label{Eq5}
    x^{5}\sim x  ~{\mathrm or}~ x^{-1}.
\end{equation}
We claim that $x^{5}\not\sim x$ for any $1 \neq x\in G$. For, if $x^{5}\sim x$, i.e., $x^{5}=g^{-1}xg$, for some $g \in G$, then $x^{5^{o(g)}}=x$, which implies \begin{equation}\label{Eq1}
    5^{o(g)}\equiv 1 ~(\mathrm{mod}~o(x)).
\end{equation}

\noindent Observe that in either of the two possibilities given by (\ref{Eq4}), for the value of $o(x)$, 5 is a primitive root of $o(x)$. Therefore, (\ref{Eq1}) yields that $\varphi(o(x))$ divides $o(g)$, where $\varphi$ denotes the Euler's phi function. This is clearly not possible as $o(g)$ is odd. It thus follows that
 \begin{equation}\label{Eq6}
    x^{5}\sim x^{- 1}~\mathrm{for~ all}~ x\in G,
\end{equation}
 i.e., (i) holds. Arguing as above, it is easy to see that if $o(x)=7^{a}$, then (\ref{Eq6}) is possible only if $a\leq 1$. Hence, (ii) holds.

\para Conversely, let $G$ be a finite group satisfying (i) and (ii). Since $5$ is a primitive root of $7$ and of $3^{a}$, for all $a\geq 1$, (i) implies that for every $x\in G$, $x^{j}\sim x$ or $x^{-1}$, for all $j$ relatively prime to $o(x)$. Hence, in view of (\ref{Eq4a}), $G$ has the {\sf cut}-property. ~$\Box$\\

 Let $V(\mathbb{Z}[G])$ denote the group of units of augmentation 1 in $\mathbb{Z}[G]$. In \cite{BMP}, a classification of  metacyclic groups has been given according to the central height of $V(\mathbb{Z}[G])$, i.e., the least integer $n$ at which the upper central series $\{\mathcal Z_i(V(\mathbb Z[G]))\}_{i\geq 0}$ of $V(\mathbb{Z}[G])$ stabilizes. In view of (\cite{saty}, Theorem 3.7) and Theorem \ref{T5}, we conclude the following:

\begin{cor}\label{c2} Let $G$ be an odd-order group. The central height of $V(\mathbb{Z}[G])$ is $1$, except when $G$ is a group with the {\sf cut}-property such that $G$ contains an element of order 7, and in this case, the central height of $V(\mathbb{Z}[G])$ is $0$.
\end{cor}

A simplifying feature for the classification of odd-order groups with the \linebreak{\sf cut}-property is that every element of such a group necessarily has prime power order. In contrast with the odd-order groups with the {\sf cut}-property, an even-order solvable group with the {\sf cut}-property does not necessarily have all its elements of prime power order. For example, the metacyclic group $$\langle a,~b~|~ a^{12}=1,~b^{2}=1,~b^{-1}ab = a^ {5}\rangle$$ has the {\sf cut}-property (\cite{BMP}, Theorem 5), while it has elements of mixed order. For the classification of even-order solvable groups with the {\sf cut}-property, we restrict to the class of solvable groups in which every element has prime power order \cite{Higman}.\\

For a finite group $G$, let $\pi(G)$ denote the set of primes dividing the order of $G$.
\begin{theorem}\label{T7}Let $G$ be a finite solvable group such that every element of $G$ has prime power order. Then, $G$ has the {\sf cut}-property if, and only if, every element $x\in G$ satisfies one of the following conditions:

 \begin{description}
  \item[~(i)~] $o(x)=2^{a}$, $a\geq 0 $ and $x^{3}\sim x$ or $ x^{-1}$;
\item[\,(ii)\,] $o(x)= 7$ or $3^{b}$, $b\geq 1$ and $x^{5}\sim x^{-1}$;
\item[(iii)] $o(x)=5$ and $x^{3}\sim x^{-1}$.
\end{description}
\end{theorem}\textbf{Proof.} Let $G$ be a finite solvable group with the {\sf cut}-property such that every element of $G$ has prime power order. Then, $G$ satisfies (\ref{Eq4a}) and therefore, by (\cite{CD}, Theorem~2), we have that $\pi(G)\subseteq\{2,3,5,7,13\}$. Let $1\neq x\in G$. If $o(x)$ is even, then clearly, $x$ satisfies (i). If $o(x)$ is odd, then $o(x)=p^{\alpha}$, where $p=3,5,7$ or $13$ and $\alpha\geq 1$, so that $o(x)$ has a primitive root, say $r$. By (\ref{Eq4a}), we have that $x^{r}\sim x^{\varepsilon}$, $\varepsilon=\pm 1$, i.e., $x^{r}=g^{-1}x^{\varepsilon}g$, for some $g\in G$ and hence $(r\varepsilon)^{o(g)}\equiv 1~(\mathrm{mod}~o(x))$. Consequently,
\newpage
\begin{equation}\label{Eq2}
    \varphi(o(x)) ~\mathrm{ divides} ~o(g), ~\mathrm{if}~ \varepsilon = 1.
\end{equation} and \begin{equation}\label{Eq3}
    \varphi(o(x)) ~\mathrm{ divides} ~2o(g), ~\mathrm{if}~ \varepsilon = -1.
\end{equation} Clearly, neither (\ref{Eq2}) nor (\ref{Eq3}) holds, if $p=13$. Thus, $G$ has no element of order 13. Moroever, if $p=7$, then (\ref{Eq2}) fails and (\ref{Eq3}) holds, only if $o(x)=7$. Furthermore, if $p=5$, then either of (\ref{Eq2}) and (\ref{Eq3}) holds, only if $\alpha=1$. Note further, that $x^{r}\sim x$ yields $x\sim x^{-1}$, if $o(x)=5$. Similarly, if $p=3$, then ($\ref{Eq2}$) holds only if $\alpha=1$ and in that case, $x^{5}\sim(=)x^{-1}$. These observations put together yield that $x$ must satisfy (ii) or (iii).

\para Converse follows easily from (\ref{Eq4a}) by observing that $5$ is a primitive root of 7 and of $3^{b}$, for all $b\geq 1$, $3$ is a primitive root of $5$, and the fact that $\mathcal{U}(\mathbb{Z}/2^{n}\mathbb{Z})=\langle \pm 1\rangle \oplus \langle 3 \rangle$, $n\geq 3$.  $\Box$\\

Recall that an element $x \in G$ is said to be real if $x\sim x^{-1}$, and $G$ is called a \emph{real group} if every element of $G$ is real. It may be of interest to note that the class of real groups with the {\sf cut}-property is the same as (i) the class of \emph{rational groups} studied by Chillag and Dolfi \cite{CD}, and (ii) the class $\mathcal{B}$ of groups studied by Golomb and Hales \cite{Hales}. By arguments given in the proof of Theorem \ref{T7}, we observe that $$\pi(G)\subseteq\{2,3,5\}, ~if~ G~ is~ a ~rational ~ solvable~ group.$$

    \section{Nilpotent groups with the {\sf cut}-property}
    If $G$ is a finite nilpotent group with the {\sf cut}-property, then $\pi(G)\subseteq \{2,3\}$; further, $2$-groups and $3$-groups with the the {\sf cut}-property have also been characterized \cite{BMP}. The following theorem extends the classification of finite abelian groups with the {\sf cut}-property \cite{Hig} to the classification of finite nilpotent groups with the {\sf cut}-property.

    \begin{theorem}\label{T4}A finite nilpotent group $G$ has the {\sf cut}-property if, and only if, $G$ is one of the following:
\begin{description}
  \item[~(i)~]a $2$-group such that for all $x \in G$, $x^{3}\sim x$ or $x^{- 1}$;
  \item[\,(ii)\,]a $3$-group such that for all $x \in G$, $x^{2}\sim x^{-1}$;
      \item[(iii)] a direct sum $H\oplus K$ of a real group $H$ satisfying~{\rm (i)} and a non-trivial group $K$ satisfying~{\rm(ii)}.

\end{description}
\end{theorem}\textbf{Proof.} Let $G$ be a finite nilpotent group with the {\sf cut}-property, so that $\pi(G)\subseteq\{2,3\}$. If $\pi(G)\neq\{2,3\}$, then $G$ is of type (i) or (ii) (\cite{BMP}, Theorems 2 \& 3). In case  $\pi(G)=\{2,3\}$, then $G=H\oplus K$, where $H$ is a non-trivial $2$-group and $K$ is a non-trivial $3$-group. Since the {\sf cut}-property is quotient closed, both $H$ and $K$ have the {\sf cut}-property and hence, $H$ satisfies (i) and $K$ satisfies (ii). It only remains to check that $H$ is a real group. For this, let order of $H$ be $2^{\alpha}$ and that of $K$ be $3^{\beta}$ for some positive integers $\alpha,~\beta$. Further, let $m$ and $n$ be positive integers satisfying $$m\equiv 3 ~(\mathrm {mod}~2^{\alpha}); ~m\equiv 1 ~(\mathrm {mod}~3^{\beta})$$ and $$n\equiv 3 ~(\mathrm {mod}~2^{\alpha});~n\equiv -1~ (\mathrm {mod}~3^{\beta}).$$ Note that neither $2$ nor $3$ can divide either $m$ or $n$ and therefore, both $m$ and $n$ are relatively prime to the order of $(h,k)$, for any $h\in H$ and $1\neq k\in K$. Since $H\oplus K$ has the {\sf cut}-property, we obtain, by (\ref{Eq4a}), that $(h,k)^{m}=(h^{3},k)$ must be conjugate to $(h,k)$ or $(h,k)^{-1}$. Now, $K$ being a $3$-group, $k\in K$ is not real. Therefore, we must have $(h,k)^{m}\sim (h,k)$ and hence $h^{3}\sim h$. Similarly, $(h,k)^{n}=(h^{3},k^{-1})$ must be conjugate to $(h,k)^{-1}$, which yields that $h^{3}\sim h^{-1}$. Consequently, $h\sim h^{-1}$, for any $h\in H$, i.e., $H$ is a real group.

\para Conversely, if $G$ is of type (i) or (ii), then clearly $G$ has the {\sf cut}-property \linebreak(\cite{BMP}, Theorems 2 \& 3). Let $G=H\oplus K$ be as in (iii) and let $(h,k)\in H \oplus K$. We check that for any integer $i$, relatively prime to the order of $(h,\,k)$, $(h,\,k)^i\sim (h,\,k)^{\pm 1}$. Since $H$ satisfies (i) and $K$ satisfies (ii), both $H$ and $K$ have the\linebreak {\sf cut}-property (\cite{BMP}, Theorems 2 \& 3) and thus satisfy (\ref{Eq4a}). Therefore, if $h=1$, or $k=1$, then $(h,\,k)^i\sim (h,\,k)^{\pm 1}$. In case neither $h=1$ nor $k=1$, then $i$ is neither divisible by 2 nor by 3. Therefore, $h^i\sim h\sim h^{-1}$ and $k^i\sim k^{\pm 1}$. Hence $(h,\, k)^{i}\sim (h,\,k)^{\pm 1}$. It now follows from (\ref{Eq4a}), that $H\oplus K$ has the {\sf cut}-property.  $\Box$\\

{\bf Remark.} It is known that the direct sum $H\oplus K$ of two $2$-groups with the {\sf cut}-property may not have the {\sf cut}-property (\cite{BMP}, Remark 1). However, following the arguments of Theorem \ref{T4}, we see that this is so if one of $H$ or $K$ is a real group. By considering the direct sum of cyclic groups of order $4$, one can easily check that the above condition is sufficient, but not necessary. More precisely, the direct sum $H\oplus K$ of 2-groups $H$ and $K$ with the {\sf cut}-property, does not have the {\sf cut}-property if, and only if, there exist non-real elements $h \in H $ and $k \in K$, such that either $h^{3}\sim h$ and $k^{3}\sim k^{-1}$ or $h^{3}\sim h^{-1}$ and $k^{3}\sim k$.\\

  In view of the proof of Theorem \ref{T4} and the above remark, we obtain the following result analogous to (\cite{Hig}, Theorem 7).

  \begin{cor}\label{P6} Let $G$ be a finite nilpotent group with the {\sf cut}-property and let \linebreak $G^{*}:=G \oplus H$, where $H$ is a real $2$-group with the {\sf cut}-property. Then, $$ \mathcal{Z}(\mathcal{U}(\mathbb{Z}[G^{*}]))=\pm\mathcal{Z}(G^{*}).$$

  \end{cor}

\para We next give two alternative characterizations for a $p$-group of class 2 to have the {\sf cut}-property.\\

For an element $x$ in a group $G$, define $[x,G]:=\{[x,g]:=xgx^{-1}g^{-1}~|~g\in G\}$. Note that in case $G$ is of class $2$, then $[x, G]$ is a normal subgroup of $G$. Thus, as a consequence of Theorem \ref{T4}, we have

 \begin{cor}\label{T1}A $p$-group $G$ of class 2 has the {\sf cut}-property if, and only if, one of the following holds:
 \begin{description}
   \item[\,(i)\,]  $p=2$ and, for every $x\in G$, $x^{4}\in [x,G]$;
        \item[(ii)] $p=3$ and, for every $x\in G$, $x^{3}\in [x,G]$.

 \end{description}
Furthermore, the {\sf cut}-property is direct sum closed for $p$-groups of class 2.
 \end{cor}
We also have the following:
 
 \begin{prop}\label{P4} A finite $p$-group $G$ of class 2 has the {\sf cut}-property if, and only if,
  \begin{equation}\label{Eq7}
  for~ all~ x\in G,~ both ~[x,G] ~and~ G/[x,G] ~ have~ the~ {\sf cut}{\rm-}property.
  \end{equation}

\end{prop}\textbf{Proof.} The necessity follows from the observation that $[x,G]\subseteq \mathcal{Z}(G)$ and the fact that the {\sf cut}-property is centre closed and quotient closed.

\para Let $G$ be a $p$-group of class 2 satisfying (\ref{Eq7}). Since $[x,G]$ is a $p$-group which has the {\sf cut}-property, $p=2$ or $3$.

\para Let $p=2$. If $x\not\in [x,G]$, then $\overline{x}:=x[x,G]\in G/[x,G]=:\overline{G}$, is such that $\overline{x}^{4}\in [\overline{x},\overline{G}]$, as $\overline{G}$ has the {\sf cut}-property. This yields that $x^{4}\in [x,G]$ and hence $G$ has the {\sf cut}-property. \para The case when $p=3$ is similar and we omit the details. $\Box$\\

 If a group $G$ has the {\sf cut}-property, then so do $\mathcal{Z}(G)$ and $G/\mathcal{Z}(G)$. The converse, however, is not true, even for the nilpotent groups of class 2. For example, if $$G=\langle a,~b~|~a^{9}=1,~b^{9}=1,~b^{-1}ab=a^{4} \rangle,$$ then both $\mathcal{Z}(G)$ and $G/\mathcal{Z}(G)$, being abelian of exponent 3, have the {\sf cut}-property, whereas $G$ does not have the {\sf cut}-property, since $b^{2}\not\sim b^{-1}$. However, Proposition \ref{P4} yields the following:
 \begin{cor}\label{T3} A finite $p$-group $G$ of class 2 has the {\sf cut}-property if, and only if, for all normal subgroups $N$ of $G$ contained in $\mathcal{Z}(G)$, both $N$ and $G/N$ have the {\sf cut}-property.
                                                      \end{cor}

\section{Concluding remark}
 Our analysis leaves open the problem of characterizing a solvable group $G$ with the {\sf cut}-property in case $G$ is a non-nilpotent solvable group of even order having an element of mixed order.\\

\noindent\textbf{Acknowledgement}\\
\noindent The author is grateful to I.\,B.\,S. Passi for his valuable comments and \linebreak suggestions.

\bibliographystyle{amsplain}
\bibliography{ReferencesM}
\end{document}